\documentclass[preprint,12pt]{elsarticle}
\newtheorem{theorem}{Theorem}
\newtheorem{corollary}{Corollary}

\usepackage{graphicx}
\usepackage{graphics}
\usepackage{subfigure}
\usepackage{epstopdf}
\usepackage{harpoon}
\usepackage{amssymb}
\usepackage{float}
\usepackage{caption2}
\usepackage{geometry}
\usepackage{picins}
\usepackage{subfigure}
\usepackage{extarrows}

\usepackage{amssymb}

\begin{document}

\begin{frontmatter}



\title{Singular cycles connecting saddle periodic orbit and saddle equilibrium in piecewise smooth systems }

\setcounter{footnote}{0}
\author{Lei Wang$^{1}$ and Xiao-Song Yang$^{2}$}

\address{$^1$Department of Mathematics and Physics, Hefei University, Hefei 230601, China\\
 $^2$School of Mathematics and Statistics, Huazhong University of Science and Technology, Wuhan 430074, China\\
$^{*}$\textup{Author for correspondence. Email: xsyang@hust.edu.cn}.}

\begin{abstract}
For flows, the singular cycles connecting saddle periodic orbit and saddle equilibrium can potentially result in the so-called singular horseshoe, which means the existence of a
non-uniformly hyperbolic chaotic invariant set. However, it is very hard to find a specific dynamical system that exhibits such singular cycles in general.
In this paper, the existence of the singular cycles involved in saddle periodic orbits is studied by two types of piecewise
affine systems: one is the piecewise
affine system having an admissible saddle point with only real eigenvalues
and an admissible saddle periodic orbit, and the other is the piecewise
affine system having an admissible saddle-focus
and an admissible saddle periodic orbit. Precisely, several kinds of sufficient conditions are obtained for the existence of only one heteroclinic cycle or only two heteroclinic
cycles in the two types of piecewise
affine systems, respectively.  In addition, some examples are
presented to illustrate the results.

\end{abstract}


\end{frontmatter}


\section{Introduction}
Singular cycle (generally refers to a homoclinic orbit connecting a singularity to itself, or a heteroclinic cycle connecting different singular elements at least one of which is a singularity\cite{Ara2010Three}) is one of the most important mechanisms leading to complicated dynamic behaviors and bifurcations \cite{Shilnikov2001Methods,Wiggins1990Introduction,Homburg2010Homoclinic,Liang2012Limit,Wei2015Limit}. For example, the well known Shil'nikov theory \cite{Shilnikov2001Methods,Wiggins1990Introduction} shows that the existence of singular cycles involved only in equilibrium points
implies the existence of a countable number of chaotic invariant sets under some
conditions. In addition, another type of singular cycle connecting a saddle equilibrium point and a saddle periodic orbit can generate a so-called singular horseshoe, which means the existence of a non-uniformly hyperbolic chaotic invariant set, and is the basis for the theoretical study of the geometrical Lorenz model\cite{Ara2010Three,Guckenheimer1979Structural,Carrasco2010One,Labarca1986Stability} motivated from the study of the famous Lorenz attractor \cite{Lorenz2004Deterministic,Afra1977Origin,Tucker1999The}.\\
\indent Although theoretical importance of singular cycles, it is very difficult to prove the existence of the singular cycles for a concrete smooth system, which is always a meaningful topic for many years \cite{BO1994CONSTRUCTING, Deng2002Food,Tigan2011Analytical,Leonov2012General,Leonov2014Fishing}. Consequently, it is interesting to find some concrete systems possessing the singular cycles.\\
\indent Due to the wide applications, the dynamics in piecewise smooth have received much attention recently \cite{leine2013dynamics,bernardo2008piecewise,llibre2003existence,llibre2013existence,wang2018existence1,wang2018existence}. By using the piecewise smooth models, the existence of some types of singular cycles involved only in equilibrium points were proved by studying the spatial location relationship between the (un)stable manifolds of the subsystems and the switching manifold, based on which some important results of the Shil'nikov theory have been generalized and upgraded to the piecewise smooth systems, see \cite{Huan2012Chaos,Wang2017Heteroclinic,Wu2016Chaos,Song2014Existence}. \\
\indent
In this paper, we will study the existence of singular cycles involved in saddle periodic orbits in a class of piecewise smooth systems. In more detail, two types of such singular cycles are investigated: one is the heteroclinic cycle connecting a saddle point with only real eigenvalues  and a saddle periodic orbit, and the other is the heteroclinic cycle connecting a saddle-focus point and a saddle periodic orbit. In particular, for each type of singular cycles, we obtain the sufficient conditions under which there exists only one heteroclinic cycle or only two heteroclinic cycles. At last, some examples are given to illustrate our main results.\\
\indent The rest of this paper is organized as follows. Section 2 introduces the considered systems and
two main results of this paper. Section 3 gives some lemmas. Sections 4 shows in detail the proofs of the main results
in Section 2. Section 5 presents some examples to illustrate the main results.
Section 6 gives some further conclusions.

\section{Systems and Main results: Existence of heteroclinic cycles connecting a periodic orbit and an equilibrium}
Consider the following systems:
 \begin{eqnarray}\label{M1}
\dot{\mathbf{x}}=\left\{\begin{array}{l}
(A-\textup{diag}(x_1^2+x_2^2,x_1^2+x_2^2,0))\mathbf{x}, ~~ \textup{if}~~\textbf{c}^T\mathbf{x}\leq d\\
B(\mathbf{x}-\textbf{q}),~~~~~~~\textup{if}~~ \textbf{c}^T\mathbf{x}>d
\end{array}\right.,
\end{eqnarray}
where $\mathbf{x}=(x_1, x_2,x_3)^T\in \mathbb{R}^3$ is a vector of state variables,
$d\in \mathbb{R}^+$ is a constant$, \textbf{q}=(q_1,q_2,q_3)^T\in\mathbb{R}^3$, $\textbf{c}=(1,0,1)^T\in\mathbb{R}^3$, and
 \begin{eqnarray*}
 A=\left(
 \begin{array}{ccc}
    \rho &-\omega &0\\
    \omega & \rho& 0\\
    0 & 0& \mu
  \end{array}\right),~B=\left(
 \begin{array}{ccc}
    b_{11} &b_{12} &0\\
    b_{21} & b_{22}& 0\\
    0 & 0& \lambda
  \end{array}\right)
  \end{eqnarray*}
with $\rho\in \mathbb{R}^+$, $\omega\in \mathbb{R}^+$, $\mu\in \mathbb{R}^+$, $\lambda\in \mathbb{R}^+$ and $b_{ij}\in \mathbb{R}(i,j=1,2)$. Moreover, matrix $B$ satisfies the following hypotheses\\
$~~~~~$\emph{(H1) the eigenvalues of $B$: $\lambda_{1,2}<0$ and $\lambda>0$,}\\
or \emph{(H2) the eigenvalues of $B$: $\alpha\pm \beta i$ and $\lambda$ where  $\alpha<0$, $\beta>0$ and $\lambda>0$};
\\ and matrix $A$ and $\mathbf{q}$ satisfy the following hypothesis \\
$~~~~~$\emph{(H3) $0<\sqrt{\rho}<d,~\textbf{c}^T\mathbf{q}>d, ~q_1=d.$}
\\
\indent For convenience, we give some notations in the following. Denote by
 \begin{eqnarray*}\label{subflow}
\phi_A(t,\cdot) ~~\textup{and}~~\phi_B(t,\cdot)
\end{eqnarray*}
 the flows generated by the left subsystem
 \begin{eqnarray}\label{M2}
\dot{\mathbf{x}}=
(A-\textup{diag}(x_1^2+x_2^2,x_1^2+x_2^2,0))\mathbf{x}, ~x\in \mathbb{R}^3
\end{eqnarray}
and right subsystem
 \begin{eqnarray}\label{M3}
\dot{\mathbf{x}}=
B(\mathbf{x}-\textbf{q}),~x\in \mathbb{R}^3
\end{eqnarray}
respectively. Obviously, $\mathbf{q}$ is the only saddle equilibrium of (\ref{M3}) with
its stable manifold and unstable manifold being
 \begin{eqnarray}\label{smq}
 W^s(\mathbf{q})=\{\mathbf{x}\in \mathbb{R}^3|x_3=q_3\}, W^u(\mathbf{q})=\{\mathbf{x}\in \mathbb{R}^3|x_1=q_1,x_2=q_2\},
 \end{eqnarray}
respectively. Furthermore, by the general polar coordinates transformation to $x_1$ and $x_2$, i.e.,
\begin{eqnarray}\label{vand}
x_1=rcos\theta,x_2=rsin\theta,
\end{eqnarray}
(\ref{M2}) can be transformed to
\begin{eqnarray*}
\left\{\begin{array}{l}
\dot{r}=r(\rho-r^2)\\
\dot{\theta}=\omega\\
\dot{x_3}=\mu x_3\\
\end{array}\right..
\end{eqnarray*}
Let \begin{eqnarray}\label{Upsilon}
\Upsilon=\{(x_1,x_2,x_3)\in \mathbb{R}^3|x_1^2+x_2^2=\rho,x_3=0\}.
\end{eqnarray}
From the classical results related to the Van der Pol oscillator \cite{Wiggins1990Introduction, Holmes1977Bifurcations} (refer the analysis for (\ref{S3}) below in Section 3),  it is not hard to see that $\Upsilon$ is the only saddle periodic orbit of (\ref{M2}) with its stable manifold and unstable manifold being
 \begin{eqnarray}\label{smUpsilon}
 W^s(\Upsilon)=\{0\neq\mathbf{x}\in \mathbb{R}^3|x_3=0\}, W^u(\Upsilon)=\{\mathbf{x}\in\mathbb{R}^3|x_1^2+x^2_2=\rho\},
 \end{eqnarray}
respectively.\\
\indent For simplifying further statements, let
\begin{eqnarray}
\Sigma=\{\mathbf{x}\in \mathbb{R}^3|\mathbf{c}^T\mathbf{x}=d\}, \Sigma^-=\{\mathbf{x}\in \mathbb{R}^3|\mathbf{c}^T\mathbf{x}<d\}, \Sigma^+=\{\mathbf{x}\in \mathbb{R}^3|\mathbf{c}^T\mathbf{x}>d\},\label{Sigma1}\\
\mathbf{p}_0=(\sqrt{\rho},0,d-\sqrt{\rho})^T\in \Sigma,~~ \mathbf{p}_1=(-\sqrt{\rho},0,d+\sqrt{\rho})^T\in \Sigma
,~~\mathbf{q}_0=(d,q_2,0)^T\in \Sigma,\label{Sigma2}\\
L_1=\{\mathbf{x}\in \mathbb{R}^3|x_1=d,x_3=0\}\subset \Sigma,~~ L_2=\{\mathbf{x}\in \mathbb{R}^3|x_1=d-q_3,x_3=q_3\}\subset \Sigma,\label{Sigma3}\\
\sigma_{\pm}=\frac{-\omega\pm\sqrt{\omega^2- 4d^2(d^2-\rho)}}{2d},\mathbf{v}_1=(d,\sigma_+,0)^T\in L_1,\label{Sigma4}\\
\mathbf{x}_-=\frac{d-\mathbf{c}^T\mathbf{q}}{\mathbf{c}^TB^{-1}\mathbf{c}^\bot} B^{-1}\mathbf{c}^\bot+\mathbf{q }~~\textup{with}~~\mathbf{c}^\bot=(0,1,0)^T.\label{Sigma5}
\end{eqnarray}
It is readily achieved that $\mathbf{c}^T\mathbf{x}_-=d$. Thus $$\mathbf{x}_-\in \Sigma.$$ 
From (H3), it is easy to see that 
\begin{center}
$\mathbf{q}\in \Sigma^+$ and $\Upsilon\subset \Sigma^+$,
\end{center}
which shows $\mathbf{q}$ and $\Upsilon$  are respectively the admissible saddle equilibrium and admissible saddle periodic orbit of (\ref{M1}).  
Moreover, denote the closed line segment and the open line segment between $\mathbf{x}_1$  and
$\mathbf{x}_2$ by $[\mathbf{x}_1,\mathbf{x}_2]$ and $(\mathbf{x}_1,\mathbf{x}_2)$ for any $\mathbf{x}_1,\mathbf{x}_2\in \mathbb{R}^3$ (or $\mathbb{R}^2$), respectively, i.e.,
\begin{eqnarray*}\label{M4}
\begin{aligned}~~[\mathbf{x}_1,\mathbf{x}_2]=\{\mathbf{x}|\mathbf{x}=\lambda\mathbf{x}_1+(1-\lambda) \mathbf{x}_2, 0\leq\lambda\leq1\},\\
(\mathbf{x}_1,\mathbf{x}_2)=\{\mathbf{x}|\mathbf{x}=\lambda\mathbf{x}_1+(1-\lambda) \mathbf{x}_2, 0<\lambda< 1\}.
\end{aligned}
\end{eqnarray*}
And let
\begin{eqnarray}\label{bkbb}
[\mathbf{x}_1,\mathbf{x}_2)=[\mathbf{x}_1,\mathbf{x}_2]-\{\mathbf{x}_2\},(\mathbf{x}_1,\mathbf{x}_2]=[\mathbf{x}_1,\mathbf{x}_2]-\{\mathbf{x}_1\}.
\end{eqnarray}

\indent Two main results on the existence of heteroclinic cycles connecting $\Upsilon$ and $\mathbf{q}$ can be presented
 in the following two theorems, which will be proved in Section 4.
\begin{theorem}
For system (\ref{M1}) with hypotheses (H1) and (H3).
\begin{enumerate}[(i)]
\item When  $d^2-\rho\geq \frac{\omega^2}{4d^2}$.
 \begin{enumerate}[a)]
\item if $q_3=d-\sqrt{\rho}$ and $\mathbf{c}^TB(\mathbf{p}_0-\mathbf{q})\geq0$, then there exists only one heteroclinic cycle connecting $\Upsilon$(see (\ref{Upsilon}))  and $\mathbf{q}$;
\item if $q_3=d+\sqrt{\rho}$ , $\omega^2\rho<\mu^2(d^2-\rho)$ and $\mathbf{c}^TB(\mathbf{p}_1-\mathbf{q})\geq0$, there exists only one heteroclinic cycle connecting $\Upsilon$  and $\mathbf{q}$;
\item if  $q_3\in (d-\sqrt{\rho}, d+\sqrt{\rho})$, $\omega^2\rho<\mu^2(d^2-\rho)$ and $\mathbf{c}^TB(\mathbf{p}_\pm-\mathbf{q})\geq0$, where $\mathbf{p}_\pm=(d-q_3,\pm\sqrt{\rho-(d-q_3)^2},q_3)$, then there exist only two heteroclinic cycles each of which connects $\Upsilon$  and  $\mathbf{q}$.
\end{enumerate}
\vspace{2mm}
 \item When $0<d^2-\rho<\frac{\omega^2}{4d^2}$.\\
 Then the flow $\phi_A(t,\mathbf{v}_1)$ will intersect with  $L_1$ under negative flight time, where $\mathbf{v}_1\in L_1$ is given in (\ref{Sigma4}). Denote by $\mathbf{v}^*=(d,v^*_2,0)^T\in L_1$  the first intersection of $\phi_A(t,\mathbf{v}_1)$ and $\mathbf{v}_1$ under negative flight time.\\
 Additionally, suppose that the following two conditions hold
\begin{enumerate}[1)]
\item if $v^*_2>\sigma_+$, then $q_2\in [\sigma_+,v^*_2]$,\item if $v^*_2<\sigma_-$, then $q_2 \in (\infty,- v^*_2]\cup [\sigma_+,+\infty)$,
  \end{enumerate}
\noindent where $\sigma_{\pm}$ is given by (\ref{Sigma4}). Then,
the three conclusions in a), b) and c) in (i) above still hold.
\end{enumerate}

\end{theorem}

\begin{theorem}
For system (\ref{M1}) with hypotheses (H2) and (H3), then $\phi_B(t,\mathbf{x}_-)$ must intersect with $L_2$ under negative flight time, where $\mathbf{x}_1$ and $L_2$ are defined as (\ref{Sigma5}) and (\ref{Sigma3}) respectively. Denote by $\mathbf{x}_+$ the first intersection of the flow $\phi_B(t,\mathbf{x}_-)$ and $L_2$ under
negative flight time.
\begin{enumerate}[(i)]
\item When  $d^2-\rho\geq \frac{\omega^2}{4d^2}$.
 \begin{enumerate}[a)]
\item if $q_3=d-\sqrt{\rho}$ and $\mathbf{p}_0\in [\mathbf{x}_-,\mathbf{x}_+)$, then there exists only one heteroclinic cycle connecting $\Upsilon$  and $\mathbf{q}$;
\item if $q_3=d+\sqrt{\rho}$, $\omega^2\rho<\mu^2(d^2-\rho)$ and $\mathbf{p}_1\in [\mathbf{x}_-,\mathbf{x}_+)$, there exists only one heteroclinic cycle connecting  $\Upsilon$ and  $\mathbf{q}$;
\item if  $q_3\in (d-\sqrt{\rho}, d+\sqrt{\rho})$, $\omega^2\rho<\mu^2(d^2-\rho)$ and $\mathbf{p}_\pm\in [\mathbf{x}_-,\mathbf{x}_+),$
where
\begin{eqnarray}\label{p+-}
\mathbf{p}_\pm=(d-q_3,\pm\sqrt{\rho-(d-q_3)^2},q_3),
\end{eqnarray}
 there exists only two heteroclinic cycle each of which connects  $\Upsilon$ and $\mathbf{q}$. Here $[\mathbf{x}_-,\mathbf{x}_+)$ is defined as (\ref{bkbb}).
\end{enumerate}

\vspace{2mm}
\item When $0<d^2-\rho<\frac{\omega^2}{4d^2}$.\\
Then the flow $\phi_A(t,\mathbf{v}_1)$ will intersect with  $L_1$ under negative flight time. Denote by $\mathbf{v}^*=(d,v^*_2,0)^T\in L_1$  the first intersection of $\phi_A(t,\mathbf{v}_1)$ and $\mathbf{v}_1$ under negative flight time.\\
 Additionally,
 suppose that the following two conditions hold
\begin{enumerate}[1)]
\item if  $v^*_2>\sigma_+$, $q_2\in [\sigma_+,v^*_2]$,
\item if $v^*_2<\sigma_-$, $q_2 \in (\infty,- v^*_2]\cup [\sigma_+,+\infty)$,
  \end{enumerate}
  then the three conclusions a), b) and c) in (i) still hold.
\end{enumerate}
\end{theorem}

\indent To accomplish the proof of Theorem 1 and Theorem 2, we need some preliminaries in Section 3.

\section{Preliminaries: some important lemmas in planar smooth systems}

\subsection{An interesting result on the classical van der Pol oscillator}
Considering the classical van der Pol oscillator
\begin{eqnarray}\label{S1}
\dot{\mathbf{x}}=(A_0-||\mathbf{x}||^2I)\mathbf{x},
\end{eqnarray}
where $\mathbf{x}=(x_1,x_2)^T\in \mathbb{R}^2$, $A_0=\left(
 \begin{array}{cc}
    \rho &-\omega\\
    \omega & \rho
  \end{array}\right)$ with $\rho>0$ and $\omega>0$, $||\mathbf{x}||=\sqrt{x_1^2+x_2^2}$ and $I$ denotes the identity matrix of order 2.\\
\indent For $\textbf{\textup{x}}_0\in \mathbb{R}^2$, denote by $O(\textbf{\textup{x}}_0)$, $O_+(\textbf{\textup{x}}_0)$ and $O_-(\textbf{\textup{x}}_0)$ the whole orbit, the positive semi-orbit and the negative semi-orbit of $\textbf{\textup{x}}_0$, respectively, i.e.,
\begin{center}
$O(\textbf{\textup{x}}_0)=\{\varphi(t,\textbf{\textup{x}}_0)|t\in \mathbb{R}\},~O_+(\textbf{\textup{x}}_0)=\{\varphi(t,\textbf{\textup{x}}_0)|t>0\}$ and $O_-(\textbf{\textup{x}}_0)=\{\varphi(t,\textbf{\textup{x}}_0)|t<0\},$
\end{center}
where $\varphi(t,\cdot)$ denotes the flow generated by (\ref{S1}).
Let $$L=\{\textbf{\textup{x}}\in \mathbb{R}^2|x_1=k\}$$ with $k>0$. Then, $L$ is perpendicular to $x_1$-axis, does not pass through the origin and divides the plane into three disjoint subsets $L$, $L_+$ and $L_-$, where
$$L_+=\{\textbf{\textup{x}}\in \mathbb{R}^2|x_1>k\},L_-=\{\textbf{\textup{x}}\in \mathbb{R}^2|x_1<k\}.$$
Obviously, the origin is in $L_{-}$. In addition, let
\begin{eqnarray}\label{varrho}
\varrho_{\pm}=\frac{-\omega\pm\sqrt{\omega^2- 4k^2(k^2-\rho)}}{2k}.
\end{eqnarray}

\begin{figure}[!htbp]
\renewcommand{\figurename}{Fig.}
 \begin{minipage}{0.25\textwidth}
\caption{An illustration for the directions of the vector fields in $L_u, L_m$ and $L_d$, which are represented by the gray short directed line segments passing through $L$. And an illustration for the the proof of subcase a) in Case (ii).}\label{fig1}
\end{minipage}
\hspace{0cm}
 \begin{minipage}{0.75\textwidth}
\includegraphics[width=14cm]{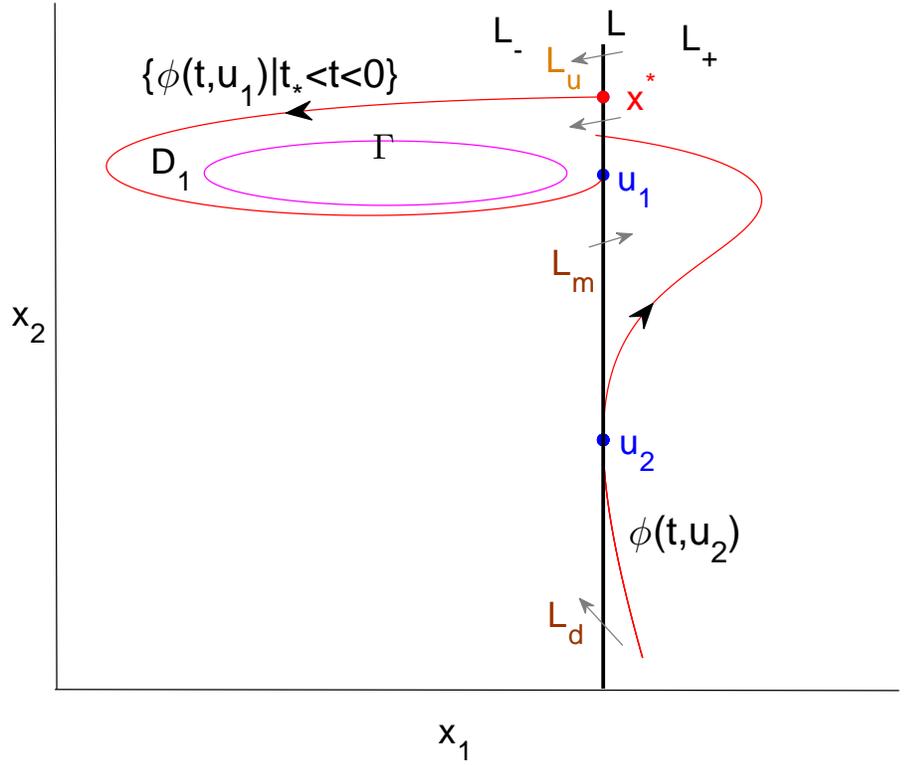}
\end{minipage}
\end{figure}

\begin{figure}[!htbp]
\renewcommand{\figurename}{Fig.}
 \begin{minipage}{0.25\textwidth}
\caption{an illustration for the the proof of subcase b) in Case (ii). Here the gray short directed line segments passing through $L$ represents the directions of vector field in $L$.}\label{fig2}
\end{minipage}
\hspace{0cm}
 \begin{minipage}{0.75\textwidth}
\includegraphics[width=12cm]{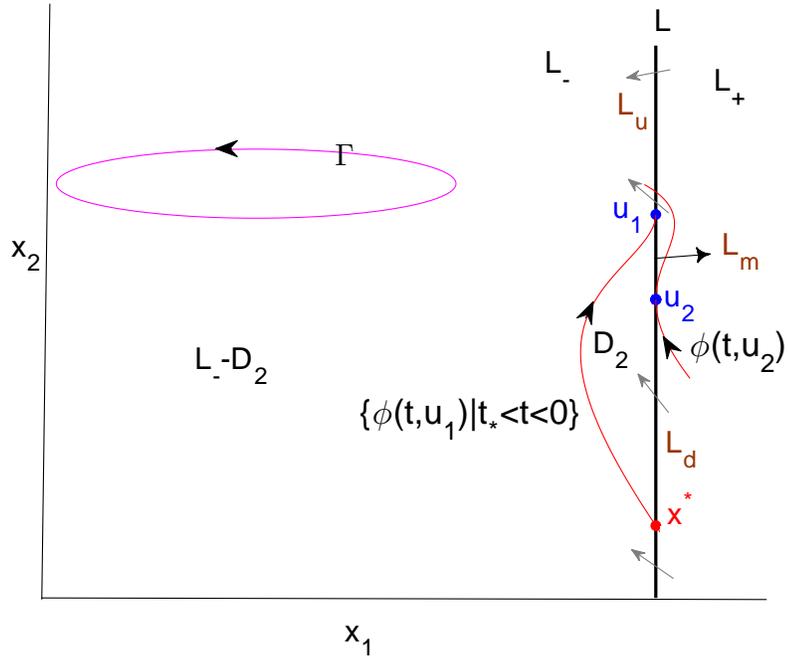}
\end{minipage}
\end{figure}

\textbf{Lemma 1}.
 \begin{enumerate}[(i)]
 \item\emph{When} $k^2-\rho\geq \frac{\omega^2}{4k^2}$.
 \begin{center}
 $O_+(\mathbf{x})\subset L_-$ \emph{for any} $\mathbf{x}\in L$.
 \end{center}
   \item \emph{When} $0<k^2-\rho<\frac{\omega^2}{4k^2}$.\\
    \emph{Let}
\begin{eqnarray*}\label{S4}
\mathbf{u}_1=(k,\varrho_+),
\mathbf{u}_2=(k,\varrho_-),
\end{eqnarray*}
then $O_-(\mathbf{u}_1)$ must intersect with $L$. \emph{Denote by} $\mathbf{x}^*=(k,x^*_2)^T\in L$ \emph{the first intersection of the flow} $\varphi(t,\mathbf{u}_1)$
 \emph{under negative flight time and} $L$.\\
 $\textup{a)}$ \emph{If }$x^*_2>\varrho_+$, \emph{then}
$$O_+(\mathbf{x})\subset L_-\cup L \Leftrightarrow \mathbf{x}\in [\mathbf{x}^*, \mathbf{u}_1],~ \textup{for}~~ \mathbf{x}\in L$$
and $$
O_+(\mathbf{x})\subset L_- \Leftrightarrow \mathbf{x}\in (\mathbf{x}^*, \mathbf{u}_1],~ \textup{for}~~ \mathbf{x}\in L.$$

$\textup{b)}$ \emph{If} $x^*_2<\varrho_-$, \emph{then}
$$O_+(\mathbf{x})\subset L_-\cup L\Leftrightarrow \mathbf{x}\in L-(\mathbf{u}_1,\mathbf{x}^*),~ \textup{for}~~ \mathbf{x}\in L$$
and
$$O_+(\mathbf{x})\subset L_-\Leftrightarrow \mathbf{x}\in L-(\mathbf{u}_1,\mathbf{x}^*],~ \textup{for}~~ \mathbf{x}\in L.$$
    \end{enumerate}
\textbf{Proof.}
In $L$, the points at which the sector field
is tangent to $L$ must meet the following equations on $\mathbf{x}$.
  \begin{eqnarray*}
  \left\{\begin{array}{l}
x_1=k\\
\rho x_1-\omega x_2-x_1(x_1^2+x_2^2)=0
\end{array}\right.,
\end{eqnarray*}
then,
  \begin{eqnarray}\label{S2}
 k x_2^2+\omega x_2+k(k^2-\rho)=0.
\end{eqnarray}
\noindent\textbf{ Case (i):} When $k^2-\rho\geq\frac{\omega^2}{4k^2}$.\\
\indent This case is trivial. In fact, since $\dot{x}_1=-kx_2^2-\omega x_2+k(\rho-k^2)\ \leq 0$ in $L$, the vector field in $L$ either points to the interior of $L_-$ or is tangent to $L$ at only one point (i.e., $(k,\frac{-\omega}{2k})^T)$. Hence, for any $\mathbf{x}\in L$, $O_+(\mathbf{x})\subset L_-$.\\\\
\noindent \textbf{Case (ii)}: When $0<k^2-\rho<\frac{\omega^2}{4k^2}$.\\
\indent In the case, $\varrho_\pm$ defined by (\ref{varrho}) are just two negative real roots of (\ref{S2}).
Thus, $\mathbf{u}_{1}$ and $\mathbf{u}_{2}$ are the only two points in $L$ at which the vector fields are tangent to $L$,
and divide $L$ into the following three segments: $ L_{u}, L_m, L_d$ (see Figs \ref{fig1}, \ref{fig2}) where
$$L_{u}=\{(x_1,x_2)^T\in L|x_2>\varrho_+\},L_{m}=\{(x_1,x_2)^T\in L|\varrho_-<x_2<\varrho_+\},L_{d}=\{(x_1,x_2)^T\in L|x_2<\varrho_-\}.$$
Obviously $L=L_{u}\cup {\mathbf{u}_1}\cup L_m\cup{\mathbf{u}_2} \cup L_d$. Now we analyse the directions of vector fields of (\ref{S1}) at $L_u,L_m,L_d$ respectively.\\
\indent 1) Provided $(x_1,x_2)^T \in L_u$.\\
 \indent Since $x_2>\varrho_+$,
we have $\dot{x}_1=-kx^2_2-\omega x_2-k(k^2-\rho)<0$ by the general nature of quadratic function, which shows that in $L_u$ the vector field of system (\ref{S1}) must  be transverse to the direction of $L$ and point to the interior of $L_-$, see Figure \ref{fig1}. \\
\indent 2) Provided $(x_1,x_2)^T \in L_m$.\\
\indent Since $\varrho_-<x_2<\varrho_+$, $\dot{x}_1=-kx^2_2-\omega x_2-k(k^2-\rho)>0$ in $L_m$. Therefore, in $L_m$ the vector field of (\ref{S1}) must be transverse to the direction of $L$ and point to the interior of $L_+$, see Figure \ref{fig1}.\\
\indent 3) Provided $(x_1,x_2)^T \in L_d$.\\
\indent Similarly, it can be readily got that in $L_u$ the vector field of (\ref{S1}) is transverse to the direction of $L$ and point to the interior of $L_-$, see Figure \ref{fig1}.\\
\indent From the discussions in 1), 2) and 3), it follows that $\varphi(t,\mathbf{u}_1)\in L_-$ and $\varphi(t,\mathbf{u}_2)\in L_+$ for small $t\neq 0$ . \\
\indent In addition, by the general polar coordinates transformation (\ref{vand}), (\ref{S1}) can be transformed to
  \begin{eqnarray}\label{S3}
  \left\{\begin{array}{l}
\dot{r}=r(\rho-r^2)\\
\dot{\theta}=\omega
\end{array}\right..
\end{eqnarray}
As we all known, (\ref{S3}) has an asymptotically stable limit cycle $\Gamma:r=\sqrt{\rho}$ with its attracting region being $\mathbb{R}^2-\{O\}$, and the flow $r(t,r_0)$ with $r_0>\sqrt{\rho}$ will tend to positive infinity  with clockwise rotation around $\Gamma$ as $t\rightarrow -\infty$.
 Thus the negative semi-orbit $O_-(\mathbf{u}_1)$ must intersect $L$ infinite times. Denote by $\mathbf{x}^*=(k,x^*_2)^T=\varphi(t_*,\mathbf{u}_1)$ the first intersection of the flow $\varphi(t,\mathbf{u}_1)$ under
negative flight time $t_*$ with $L$. Then, from the direction of vector field in $L$ showed in 1), 2) and 3), $\mathbf{x}^*$ must be either in $L_u$ (i.e., $x^*_2>\varrho_+$, see the red point in Figure 1) or in $L_d$ (i.e., $x^*_2<\varrho_-$, see the red point in Figure 2).\\
 \indent a) If $\mathbf{x}^*\in L_u$ (i.e., $x^*_2>\varrho_+$), see Figure 1.\\
\indent Consider the open region $D_1$ surrounded by $[\mathbf{x}^*, \mathbf{u}_1]\cup \{\varphi(t,\mathbf{u}_1)|t_*<t<0\}$ as shown in Figure \ref{fig1}. Obviously, $D$ is a positively invariant set contained in $L_-$. Therefore, for any $\mathbf{x} \in [\mathbf{x}^*, \mathbf{u}_1]$, $O_+(x) \subset L_-$ and $\phi(t,\mathbf{x})\rightarrow\Gamma$ as $t\rightarrow +\infty$.
In the other hand, since the attracting region of $\Gamma$ is $\mathbb{R}^2-\{O\}$, for any $\mathbf{x}\in (L_u-[\mathbf{x}^*, \mathbf{u}_1])\cup L_d$,
the flow $\varphi(t,\mathbf{x})$ must first enter into $L_--D$, and then leave $L_-$ by intersecting $L_m$,
and enter into $D_1$ ultimately by intersecting $(\mathbf{x}^*, \mathbf{u}_1)$ under positive flight time.
Moreover, it is obvious that for $\mathbf{x}\in L_m$,
 $\varphi(t,\mathbf{x})$ will first enter into $L_+$ and eventually enter into $D_1$ by
 intersecting $(\mathbf{x}^*, \mathbf{u}_1)$ under positive flight time. In conclusion,
$$O_+(\mathbf{x})\subset L_-\Leftrightarrow \mathbf{x}\in (\mathbf{x}^*, \mathbf{u}_1],~ \textup{for}~~ \mathbf{x}\in L.$$
Furthermore, $O_+(\mathbf{x}^*)\subset L_-\cup L$ is tangent to $L$ at $\mathbf{u}_1$. Thus
$$O_+(\mathbf{x})\subset L_-\cup L \Leftrightarrow \mathbf{x}\in [\mathbf{x}^*, \mathbf{u}_1],~ \textup{for}~~ \mathbf{x}\in L.$$
\indent  b) If $\mathbf{x}^*\in L_d$ (i.e., $x^*_2<\varrho_-$), see Figure 2.\\
\indent Consider the open region $D_2$ surrounded by $[\mathbf{u}_1,\mathbf{x}^*]\cup \{\varphi(t,\mathbf{u}_1)|t_*<t<0\})$. Obviously, $L_--D_2$ is a positively invariant set contained in $L_-$. Furthermore, for any $\mathbf{x}\in (\mathbf{u}_2, \mathbf{x}^*)$, under positive flight time $\varphi(t,\mathbf{x})$ must enter first into $D_2$, and then enter into $L_+$ by intersecting $L_m$, and eventually enter into $L_-$ by intersecting $L_u$. Combining the directions of the vector fields in $L$ shown in 1), 2) and 3) above, it follows that
$$O_+(\mathbf{x})\subset L_-\cup L\Leftrightarrow \mathbf{x}\in L-(\mathbf{u}_1,\mathbf{x}^*),~ \textup{for}~~ \mathbf{x}\in L$$
and
$$O_+(\mathbf{x})\subset L_-\Leftrightarrow \mathbf{x}\in L-(\mathbf{u}_1,\mathbf{x}^*],~ \textup{for}~~ \mathbf{x}\in L.$$
~~~~~~~~$\square$\\\\
\indent From the proof of Lemma 1, it is not hard to get the following conclusion if omitting the points in $L$ at which the vector fields are tangent to $L$ in all cases.
\begin{corollary}~~~~~~~~~~~~~~~~~~~~~~~~~~~~~~~~~~~~~~~~~~~~~~~~~~~~~~~~~~~~~~~~~~
\begin{enumerate}[(i)]
 \item When $0<k^2-\rho<\frac{\omega^2}{4k^2}$.   \\
 \textup{a)} If $x^*_2>\varrho_+$, then, for $\mathbf{x}\in L$, $O_+(\mathbf{x})\subset L_- $ and $O(\mathbf{x})$ intersects $L$ at $\mathbf{x}$ transversely if and only if
$$ \mathbf{x}\in (\mathbf{x}^*, \mathbf{u}_1).$$

\textup{b)} \emph{If} $x^*_2<\varrho_-$, then, for $\mathbf{x}\in L$, $O_+(\mathbf{x})\subset L_- $ and $O(\mathbf{x})$ intersects $L$ at $\mathbf{x}$ transversely if and only if
$$O_+(\mathbf{x})\subset L_-\Leftrightarrow \mathbf{x}\in L-[\mathbf{u}_1,\mathbf{x}^*].$$
  \item When $k^2-\rho= \frac{\omega^2}{4k^2}$, then, for $\mathbf{x}\in L$, $O_+(\mathbf{x})\subset L_- $ and $O(\mathbf{x})$ intersects $L$ at $\mathbf{x}$ transversely if and only if $$\mathbf{x}\in L-\{(k,\frac{-\omega}{2k})^T\}.$$
\item When $k^2-\rho>\frac{\omega^2}{4k^2}$, then, for any $\mathbf{x}\in L$, $O_+(\mathbf{x})\subset L_- $ and $O(\mathbf{x})$ intersects $L$ at $\mathbf{x}$ transversely.
 \end{enumerate}
\end{corollary}

\indent Now, we generalize the conclusion of Lemma 1 to the more general situation, namely, the
location relationship between any planar line $\tilde{L}$ and the flow of (\ref{S1}) with initial conditions being in $\tilde{L}$, where
\begin{eqnarray}\label{S5}
\tilde{L}=\{\textbf{\textup{x}}\in \mathbb{R}^2|\textbf{\textup{k}}^T\textbf{\textup{x}}=1\}, \mathbf{0}\neq\textbf{\textup{k}}=(k_1,k_2)^T\in \mathbb{R}^2.
\end{eqnarray}
It is readily derived  that $\tilde{L}$ can be transformed to the form of $L$  by an orthogonal coordinate transformation such that
the more general situation can be studied by Lemma 1 directly. \\
\indent Let
\begin{eqnarray}\label{S6}
\tilde{L}_+=\{\textbf{\textup{x}}\in \mathbb{R}^2|\textbf{\textup{k}}^T\textbf{\textup{x}}>1\},
\tilde{L}_-=\{\textbf{\textup{x}}\in \mathbb{R}^2|\textbf{\textup{k}}^T\textbf{\textup{x}}<1\}.
\end{eqnarray}
Then $\mathbb{R}^2$ is equivalent to the disjoint union of $\tilde{L}$, $\tilde{L}_+$ and $\tilde{L}_-$.
Moreover, denote
$$\tilde{k}=\frac{1}{\sqrt{k_1^2+k_2^2}},~~\tilde{\varrho}_\pm=\frac{-\omega\pm\sqrt{\omega^2-4\tilde{k}^2(\tilde{k}^2-\rho)}}{2\tilde{k}},$$
$$L_2=\{\mathbf{x}\in \mathbb{R}^2|x_1=\tilde{k}\},~~\mathbf{u}_1=(\tilde{k},\tilde{\varrho}_+)^T\in L_2,$$
 \begin{eqnarray}\label{B}
 B=\left(\begin{array}{cc}
           k_1\tilde{k} &-k_2\tilde{k}\\
           k_2\tilde{k} & k_1\tilde{k}
 \end{array}\right).\end{eqnarray}
Then, simple calculation shows that $$\tilde{\mathbf{u}}_1=B\mathbf{u}_1\in \tilde{L}.$$
\indent As a consequence, we have the following corollary
   \begin{corollary}~~~~~~~~~~~~~~~~~~~~~~~~~~~~~~~~~~~~~~~~~~~~~~~~~~~~~~~~~~~~~~~~
\begin{enumerate}[(i)]
 \item\emph{ When} $\tilde{k}^2-\rho\geq \frac{\omega^2}{4\tilde{k}^2}$.
 \begin{center}
 $O_+(\mathbf{x})\subset \tilde{L}_-$ \emph{for any} $\mathbf{x}\in \tilde{L}$.
 \end{center}
 \item When $0<\tilde{k}^2-\rho<\frac{\omega^2}{4\tilde{k}^2}$.
 \\ Let
$\tilde{\mathbf{x}}^*=B\mathbf{x}^*$, with $\mathbf{x}^*=(\tilde{k},x^*_2)^T\in L_2$ be the first intersection
          of the flow $\phi(t,\mathbf{u}_1)$ under negative flight time and $L_2$ .

\textup{a)} If $x^*_2>\tilde{\varrho}_+$, \emph{then}
 $$O_+(\mathbf{x})\subset \tilde{L}_-\Leftrightarrow \mathbf{x}\in [\tilde{\mathbf{x}}^*, \tilde{\mathbf{u}}_1],~ \textup{for}~~ \mathbf{x}\in \tilde{L}.$$
\textup{b)} If $x^*_2<\tilde{\varrho}_-$, \emph{then}
 $$O_+(\mathbf{x})\subset \tilde{L}_-\Leftrightarrow \mathbf{x}\in \tilde{L}-(\tilde{\mathbf{u}}_1,\tilde{\mathbf{x}}^*],~ \textup{for}~~ \mathbf{x}\in \tilde{L}.$$
\end{enumerate}
\end{corollary}
\textbf{Proof.} Obviously, $B$ defined as (\ref{B}) is an orthogonal matrix under which $L_2$ can be transformed  to $\tilde{L}$ with distance preserving. Then, by Lemma 1,
the proof of this corollary is trivial.~~~~~~~~~~~~$\square$

\subsection{Two useful results on planar linear systems}
 Consider general planar linear system as follows
\begin{eqnarray}\label{A0}
\dot{\textbf{\textup{x}}}=A_0\textbf{\textup{x}},~~\textbf{\textup{x}}\in \mathbb{R}^2.
\end{eqnarray}
\noindent\textbf{Lemma 2}\cite{Wang2017Heteroclinic}
 \emph{For system (\ref{A0}), suppose that  the eigenvalues of  ${A}_0$ are given by $\mu_{1,2}<0$ and $\textbf{\textup{x}}_0\in \tilde{L}$,  then
\noindent\begin{center}
$O_+(\textbf{\textup{x}}_0)\subset \tilde{L}_-$ if and only if $\textbf{\textup{k}}^TA_0\textbf{\textup{x}}_0\leq0.$
\end{center}}
\noindent \emph{Here $\tilde{L}$ and $\tilde{L}^-$ are defined as (\ref{S5}) and (\ref{S6}), respectively}. \\\\
\noindent
\textbf{Lemma 3}\cite{Wang2017Heteroclinic} \emph{ For system (\ref{A0}), suppose that the eigenvalues of  $\textup{A}_0$ are given by
\begin{center}
$\alpha\pm \beta i$
\end{center}
with  $\alpha<0$, $\beta>0$ and $i=\sqrt{-1}$.
Let
\begin{equation*}
\textbf{\textup{x}}_*=\frac{1}{\textbf{\textup{k}}^TA_0^{-1}\textbf{\textup{k}}^{\bot}}
A_0^{-1}\textbf{\textup{k}}^{\bot}
\end{equation*}
with $\textbf{\textup{k}}^{\bot}=(-k_2,k_1)^T$.
 Obviously $\textbf{\textup{x}}_*\in \tilde{L}$. Moreover, denote the first intersection of the flow with initial
condition $\textbf{\textup{x}}_*$  under negative flight time and straight line $\tilde{L}$ by $\textbf{\textup{x}}^*$.
Then, for $\textbf{\textup{x}}_0\in \tilde{L}$,
\begin{center}
$O_+(\textbf{\textup{x}}_0)\subset \tilde{L}_-$
 if and only if $\textbf{\textup{x}}_0\in [\textbf{\textup{x}}_*,\textbf{\textup{x}}^*),$
\end{center}}
\noindent where $[\mathbf{x}_*,\mathbf{x}^*)$ is defined as (\ref{bkbb}).

\section{Proof of Theorem 1,2}
\subsection{Proof of Theorem 1}
 \vspace{3mm} \textbf{The proof of (i)( When $d^2-\rho\geq \frac{\omega^2}{4d^2}$).} \vspace{3mm}\\
 To prove the existence of a heteroclinic cycle connecting $\Upsilon$ and $\mathbf{q}$, the key is to prove the existence of two heteroclinic orbits.\\
 From (\ref{smq}), (\ref{smUpsilon}), (\ref{Sigma1}) and (\ref{Sigma2}), we have
\begin{eqnarray}\label{q0}
\{ \mathbf{q}_0\}= W^s(\Upsilon)\cap \Sigma \cap W^u(\mathbf{q}).
\end{eqnarray}
 Then,
\begin{eqnarray*}\label{q01}
\phi_B(t,\mathbf{q}_0)\rightarrow \mathbf{q} ( t\rightarrow-\infty),\phi_A(t,\mathbf{q}_0)\rightarrow \Upsilon ( t\rightarrow+\infty).
\end{eqnarray*}
Let
\begin{eqnarray}\label{gama1}
\Gamma_1=\{\phi_B(t,\mathbf{q}_0)|t<0\}\cup \{\mathbf{q}_0\}\cup \{\phi_A(t,\mathbf{q}_0)|t>0\}.
\end{eqnarray}
We now want to prove that $\Gamma_1$ is the only heteroclinic orbit from $\mathbf{q}$ to $\Upsilon$ for system (\ref{M1}). To do this, it is sufficient to prove
 \begin{eqnarray}
\{\phi_B(t,\mathbf{q}_0)|-\infty<t<0\}\subset \Sigma^+\label{M7}\\
\{\phi_A(t,\mathbf{q}_0)|0<t<+\infty\}\subset \Sigma^-\cup \Sigma,\label{M8}
 \end{eqnarray}
From (\ref{q0}) and (\ref{q01}), it is obvious that $\{\phi_B(t,\mathbf{q}_0)|-\infty<t<0\}=(\mathbf{q},\mathbf{q}_0)$ must belong to $\Sigma^+$
 since $\mathbf{q}_0\in \Sigma$  and $\mathbf{q}\in \Sigma^+$. Thus (\ref{M7}) holds. In addition, we have
 \begin{center}
$\mathbf{q}_0\in L_1$ and $\{\phi_A(t,\mathbf{q}_0)|0<t<+\infty\}\subset W^s(\Upsilon)$
,\end{center}
 where $L_1$ is defined by (\ref{Sigma3}). Obviously, $L_1=W^s(\Upsilon)\cap \Sigma$.
From system (\ref{M2}), we know that in $W^s(\Upsilon)$ the flow of (\ref{M2}) is determined absolutely only by the planar Van der Pol oscillator
$$\dot{x}_1=\rho x_1-\omega x_2-x_1(x_1^2+x_2^2),\dot{x}_2=\omega x_1+\rho x_2-x_2(x_1^2+x_2^2).$$
Since $d^2-\rho\geq \frac{\omega^2}{4d^2}$, according to the conclusion (i) in Lemma 1, $\{\phi_A(t,\mathbf{x})|0<t<\infty\}$ is contained in $L_1^{-}=\{\mathbf{x}\in \mathbb{R}^3|x_1<d,x_3=0\}\subset \Sigma^{-}$  for any $\mathbf{x}\in L_1$ which shows (\ref{M8}) holds.
\noindent Thus $\Gamma_1$ given by (\ref{gama1}) is indeed a heteroclinic orbit from $\mathbf{q}$ to $\Upsilon$.\\
 \indent To show the existence of heteroclinic cycles, we now need to show the existence of other heteroclinic orbits of (\ref{M1}) from $\Upsilon$ to $\mathbf{q}$.
\\
\indent a) If $q_3=d-\sqrt{\rho}$ and $\mathbf{c}^TB(\mathbf{p}_0-\mathbf{q})\geq0$.
\\
From (\ref{smq}), (\ref{smUpsilon}), (\ref{Sigma1}) and (\ref{Sigma2}),
 \begin{eqnarray}\label{p0}
 \{\mathbf{p}_0\}= W^u(\Upsilon)\cap \Sigma \cap W^s(\mathbf{q}).
 \end{eqnarray}
  Then
\begin{eqnarray*}\label{q01}
\phi_B(t,\mathbf{p}_0)\rightarrow \Upsilon ( t\rightarrow-\infty),\phi_A(t,\mathbf{p}_0)\rightarrow \mathbf{q} ( t\rightarrow+\infty).
\end{eqnarray*}
Let
\begin{eqnarray}\label{gama2}
\Gamma_2=\{\phi_A(t,\mathbf{p}_0)|-\infty<t<0\}\cup \{\mathbf{p}_0\} \cup\{\phi_B(t,\mathbf{p}_0)|0<t<+\infty\}.
\end{eqnarray}

We now show that $\Gamma_2$ is a heteroclinic orbit from $\Upsilon$ to $\mathbf{q}$. For this, it is sufficient to show that
  \begin{eqnarray}
\{\phi_A(t,\mathbf{p}_0)|-\infty<t<0\}\subset \Sigma^-,\label{M6}\\
\{\phi_B(t,\mathbf{p}_0)|0<t<+\infty\}\subset \Sigma^+.\label{M9}
 \end{eqnarray}
From (\ref{p0}), it follows that for any $t_1<0$,
 $$\phi_A(t_1,\mathbf{p}_0)=(x_1(t_1), x_2(t_1),x_3(t_1))^T \subset W^u(\Upsilon).$$
Thus $x_1^2(t_1)+x_2^2(t_1)=\rho$. Additionally, $0<x_3(t_1)<d-\sqrt{\rho}$ since the $x_3$-coordinate of $\mathbf{p}_0$ is $d-\sqrt{\rho}$ from (\ref{Sigma2}) and $\dot{x}_3=\mu>0$ for any $(x_1(t),x_2(t),x_3(t))^T\in W^u(\Upsilon)$ from (\ref{M2}).
Hence
$$x_3(t_1)+x_1(t_1)\leq x_3(t_1)+\sqrt{x_1^2(t_1)+x_2^2(t_1)}<d-\sqrt{\rho}+\sqrt{\rho}=d,$$
which means that (\ref{M6}) holds.\\
\indent Furthermore,
since $\mathbf{p}_0\in W^s(\mathbf{q})\cap \Sigma$, $\{\phi_B(t,\mathbf{p}_0)|0<t<+\infty\}\subset W^s(\mathbf{q})$.
From ($\ref{M3}$), in $W^s(\mathbf{q})$, the flow of (\ref{M3}) is absolutely determined by the following planar system
 $$\left(\begin{array}{c}
   \dot{x}_1\\
    \dot{x}_2
  \end{array}\right)=B_0
  \left(\begin{array}{c}x_1-q_1\\x_2-q_2\end{array}\right),$$
where $B_0=\left(\begin{array}{cc}b_{11}&b_{12}\\b_{21}&b_{22}\end{array}\right)$ with eigenvalues $\lambda_{1,2}<0$ from (H1).
 Since $\mathbf{c}^TB(\mathbf{p}_0-\mathbf{q})\geq0$, by normal coordinate transformation, it is not hard to know (\ref{M9}) holds by using Lemma 2. Therefore, $\Gamma_2$ given by (\ref{gama2}) is ideed the only heteroclinic orbit from $\Upsilon$ to $\mathbf{q}$.\\
\indent Hence, in this case, there exists only one heteroclinic cycle connecting $\Upsilon$ and $\mathbf{q}$, which can be expressed as
 \begin{eqnarray}\label{heter1}
 \Upsilon \cup \Gamma_2 \cup \mathbf{q}\cup \Gamma_1.
 \end{eqnarray}\\
 \indent b) if  $q_3=d+\sqrt{\rho}$ , $\omega^2\rho<\mu^2(d^2-\rho)$ and $\mathbf{c}^TB(\mathbf{p}_1-\mathbf{q})\geq0$.\\\\
 In this subcase, $$\{\mathbf{p}_1\}= W^u(\Upsilon)\cap \Sigma \cap W^s(\mathbf{q}).$$
 Therefore,
 \begin{eqnarray*}\label{p10}
\phi_A(t,\mathbf{p}_1)\rightarrow \Upsilon ( t\rightarrow-\infty),\phi_B(t,\mathbf{p}_1)\rightarrow \mathbf{q} ( t\rightarrow+\infty).
\end{eqnarray*}
 Let \begin{eqnarray} \label{gama3}
 \Gamma_3=\{\phi_A(t,\mathbf{p}_1)|-\infty<t<0\}\cup \{\mathbf{p}_1\} \cup\{\phi_B(t,\mathbf{p}_1)|0<t<+\infty\}.
  \end{eqnarray}
To prove thta $\Gamma_3$ is the only heteroclinic orbit from $\Upsilon$ to $\mathbf{q}$,  we need only to show that
 \begin{eqnarray}
\{\phi_B(t,\mathbf{p}_1)|0<t<+\infty\}\subset \Sigma^+, \label{M10}\\
\{\phi_A(t,\mathbf{p}_1)|-\infty<t<0\}\subset \Sigma_-.\label{M11}
 \end{eqnarray}
 Since $\mathbf{p}_1\in W^s(\mathbf{q})\cap \Sigma$ and $\mathbf{c}^TB(\mathbf{p}_1-\mathbf{q})\geq0$, the proof of (\ref{M10}) can be easily carried out similarly to the proof (\ref{M9}).\\
\indent Now we prove (\ref{M11}). Let
$$E=W^u(\Upsilon)\cap \Sigma.$$
Then $E$ is an elliptic secant line of $W^u(\Upsilon)$ and can be parameterized by
$$E=\{(\sqrt{\rho}cos\varsigma,\sqrt{\rho}sin\varsigma, d-\sqrt{\rho}cos\varsigma)^T|\varsigma\in [0,2\pi)\}.$$
Furthermore, it is not hard to see that $\mathbf{p}_1$ is one of the endpoints of the large axis of $E$.
Moveover, for any $$\mathbf{p}(\varsigma)=(\sqrt{\rho}cos\varsigma,\sqrt{\rho}sin\varsigma, d-\sqrt{\rho}cos\varsigma)^T\in E,$$
the tangent sector to the elliptic secant line in $\Sigma$ can be calculated as
$$\mathbf{p}'(\varsigma)=(-\sqrt{\rho}sin\varsigma,\sqrt{\rho}cos\varsigma, \sqrt{\rho}sin\varsigma)^T.$$
Let $\theta_1(\varsigma)$ be the angle between $\mathbf{p}'(\varsigma)$ and $(0,0,-1)^T$. Then
$$cos(\theta_1(\varsigma))=\frac{-\sqrt{\rho}sin\varsigma}{\sqrt{\rho+\rho sin^2\varsigma}}=-\frac{sin\varsigma}{\sqrt{1+sin^2\varsigma}}.$$
Moreover, from (\ref{M2}), the sector field at $\mathbf{p}(\varsigma)$ can be formulated by
$$\mathbf{f}=(-\omega\sqrt{\rho}sin\varsigma,\omega\sqrt{\rho}cos\varsigma,\mu(d-\sqrt{\rho}cos\varsigma)).$$
Then
$$cos(\theta_2(\varsigma))=\frac{\mu(d-\sqrt{\rho}cos\varsigma)}{\sqrt{\omega^2\rho+\mu^2(d-\sqrt{\rho}cos\varsigma)^2}},$$
where $\theta_2(\varsigma)$ denotes the angle between $-\mathbf{f}$ and $(0,0,-1)^T$.
From $d>\sqrt{\rho}$, we have $cos(\theta_2(\varsigma))>0$ for any $\varsigma\in [0,2\pi)$. Thus,
 \begin{eqnarray}\label{theta2}
 \theta_2(\varsigma)\in (0,\frac{\pi}{2}).
 \end{eqnarray}
Now, we consider the sign of
$cos^2(\theta_2(\varsigma))-cos^2(\theta_1(\varsigma))$. Simple calculation shows that
 \begin{eqnarray}\label{M12}
 \begin{aligned} cos^2(\theta_2(\varsigma))-cos^2(\theta_1(\varsigma))&=\frac{\mu^2(d-\sqrt{\rho}cos\varsigma)^2-\omega^2\rho sin^2\varsigma}{(1+sin^2\varsigma)(\omega^2\rho+\mu^2(d-\sqrt{\rho}cos\varsigma)^2)}\\
 &=\frac{(\mu d-\mu\sqrt{\rho}cos\varsigma-\omega\sqrt{\rho}sin\varsigma)(\mu d-\mu\sqrt{\rho}cos\varsigma+\omega\sqrt{\rho}sin\varsigma)}{(1+sin^2\varsigma)(\omega^2\rho+\mu^2(d-\sqrt{\rho}cos\varsigma)^2)}.
 \end{aligned}
 \end{eqnarray}
 Since $\omega^2\rho<\mu^2(d^2-\rho)$, we have $\sqrt{\mu^2+\omega^2}<\frac{\mu d}{\sqrt{\rho}}$. Thus
\begin{center}
 $\mu\sqrt{\rho}cos\varsigma+\omega\sqrt{\rho}sin\varsigma\leq \sqrt{\rho}\sqrt{\mu^2+\omega^2}<\mu d$ and $\mu\sqrt{\rho}cos\varsigma-\omega\sqrt{\rho}sin\varsigma\leq \sqrt{\rho}\sqrt{\mu^2+\omega^2}<\mu d.$
\end{center}
Combining these with (\ref{M12}) show $cos^2(\theta_2(\varsigma))-cos^2(\theta_1(\varsigma))>0$, which implies that
\begin{center}
$\theta_2(\varsigma)< \theta_1(\varsigma)$
\end{center}
by (\ref{theta2}). This shows that the flow with any initial condition $\mathbf{p}(\varsigma)$ will tend to
  $\Upsilon$ along $W^s(\Upsilon)$ without intersecting with $\Sigma$ once again as $t\rightarrow-\infty$. Hence (\ref{M11}) holds.
Obviously, the only heteroclinic cycle connecting $\Upsilon$ and $\mathbf{q}$ in this subcase can be expressed as
 \begin{eqnarray}\label{heter2}
 \Upsilon\cup\Gamma_3\cup \{\mathbf{q}\}\cup\Gamma_1. \end{eqnarray}
\indent c) If $d-\sqrt{\rho}<q_3<d+\sqrt{\rho}$, $\omega^2\rho<\mu^2(d^2-\rho)$ and $\mathbf{c}^TB(\mathbf{p}_\pm-\mathbf{q})\geq0$. \\\\
In this case, it is easy to get
 \begin{eqnarray}\label{mani1}
 W^u(\Upsilon)\cap \Sigma \cap W^s(\mathbf{q})=\{\mathbf{p}_+,\mathbf{p}_-\}.
  \end{eqnarray}
from (\ref{smq}), (\ref{smUpsilon}), (\ref{Sigma1}) and (\ref{p+-}). Hence
 \begin{eqnarray}\label{q01}
\phi_A(t,\mathbf{p}_{\pm})\rightarrow \Upsilon ( t\rightarrow-\infty),\phi_B(t,\mathbf{p}_{\pm})\rightarrow \mathbf{q} ( t\rightarrow+\infty).
\end{eqnarray}
Let
 \begin{eqnarray}\label{Gammapn}
 \Gamma_{\pm}=\{\phi_A(t,\mathbf{p}_\pm)|-\infty<t<0\}\cup \{\mathbf{p}_{\pm}\} \cup\{\phi_B(t,\mathbf{p}_\pm)|0< t<+\infty\}.
 \end{eqnarray}
Since $\omega^2\rho<\mu^2(d^2-\rho)$, by the similar proof for (\ref{M11}) in b) above, we obtain
 \begin{eqnarray}\label{phi1}
 \{\phi_A(t,\mathbf{p}_\pm)|-\infty<t<0\}\subset \Sigma_-.
 \end{eqnarray}
  Since $\mathbf{c}^TB(\mathbf{p}_\pm-\mathbf{q})\geq0$,
\begin{eqnarray}\label{phi2}
\{\phi_B(t,\mathbf{p}_\pm)|0<t<+\infty\}\subset \Sigma^+
 \end{eqnarray}
 by Lemma 2.
 Then $\Gamma_{\pm}$ are two different heteroclinic orbits from $\Upsilon$ to $\mathbf{q}$ from (\ref{mani1})$\sim$(\ref{phi2}).
  Then,
 \begin{eqnarray}\label{heter3}
 \Upsilon\cup\Gamma_{\pm} \cup \{\mathbf{q}\}\cup
\Gamma_1.
\end{eqnarray} are only two heteroclinic cycles of (\ref{M1}), each of which connect $\Upsilon$ and $\mathbf{q}$.\\
\indent \vspace{3mm} The proof of (i) is accomplished.\\
\indent \textbf{The proof of (ii) (When  $0<d^2-\rho< \frac{\omega^2}{4d^2}$).} \vspace{3mm} \\
\indent Now, we first prove that $\Gamma_1$ defined still as (\ref{gama1}) is a heteroclinic orbit from $\mathbf{q}$ to $\Upsilon$ for system (\ref{M1}) in the following. To do this, it is still sufficient to prove that (\ref{M7}) and (\ref{M8}) are both hold. (\ref{M7}) is obvious ture by using the same disccusion in (i). Now it is crucial to prove (\ref{M8}). Since $0<d^2-\rho< \frac{\omega^2}{4d^2}$, $\{\phi_A(t,\mathbf{v}_1)|t<0\}$ will intersect with $L_1$ according to Lemma 1, where $\mathbf{v}_1$ is defined in (\ref{Sigma4}). Let $\mathbf{v}^*=(d,v^*_{2},0)^T$ be the first intersection of flow $\phi_A(t,\mathbf{v}_1)$ under negative flight time and $L_1$.
\\
\indent 1) if $v^*_{2}>\sigma_+$,\\
then $q_2\in [\sigma_+,v_2^*]$ which means $\mathbf{q}_0\in [\mathbf{v}_*,\mathbf{v}_1]$. This shows that (\ref{M8}) holds by Lemma 1,\\
\indent 2) if $v^*_{2}<\sigma_-$,\\
then $q_2 \in (\infty,- v^*_2]\cup [\sigma_+,+\infty)$ which means $\mathbf{q}_0\in L_1-(\mathbf{v}_1,\mathbf{v}_*)$. Hence (\ref{M8}) still holds from Lemma 1. \\
\indent Thus $\Gamma_1$ is still a heteroclinic orbit from $\Upsilon$ to $\mathbf{q}$ in this case. In addition,\\
\indent a) If $q_3=d-\sqrt{\rho}$ and $\mathbf{c}^TB(\mathbf{p}_0-\mathbf{q})\geq0$.\\
By the completely same discussion used in subcase a) in the proof of (i), $\Gamma_2$ given in (\ref{gama2}) is still a heteroclinic orbit from $\mathbf{q}$ to $\Upsilon$. Hence,
(\ref{heter1}) is still the only heteroclinic cycle of (\ref{M1}) connecting $\Upsilon$ and $\mathbf{q}$.\\
 \indent b) if  $q_3=d+\sqrt{\rho}$ , $\omega^2\rho<\mu^2(d^2-\rho)$ and $\mathbf{c}^TB(\mathbf{p}_1-\mathbf{q})\geq0.$\\
 The discussion is also completely same as subcase b) in the proof of (i). (\ref{heter2}) is the only heteroclinic cycle of (\ref{M1}) connecting $\Upsilon$ and $\mathbf{q}$ as before.\\
 \indent c) If $d-\sqrt{\rho}<q_3<d+\sqrt{\rho}$, $\omega^2\rho<\mu^2(d^2-\rho)$ and $\mathbf{c}^TB(\mathbf{p}_\pm-\mathbf{q})\geq0$ . \\
Using the same discussion as subcase c) in the proof of (i), (\ref{heter3}) are the only two heteroclinic cycles of (\ref{M1}), each of which connects $\Upsilon$ and $\mathbf{q}$.\\
 \indent The proof of (ii) is accomplished~~~~~~~~~~$\Box$.
\subsection{Proof of Theorem 2}
Under condition (H2), the fundament difference between Theorem 1 and Theorem 2 is that $\mathbf{q}$ is a pure saddle point (with purely real
eigenvalues) in Theorem 1, while $\mathbf{q}$ is a saddle-focus in Theorem 2. In addition,
Lemma 3 can play the similar role as Lemma 2 in the proof of Theorem 1 when proving Theorem 2. Thus, for the proof of Theorem 2, we don't present the more details but a briefly description below.\vspace{3mm}\\
\textbf{The proof of (i)( When $d^2-\rho\geq \frac{\omega^2}{4d^2}$).}\vspace{3mm}\\
All of the dicussions for (\ref{q0}), (\ref{gama1}), (\ref{M7}) and (\ref{M8}) can be achieved in this case. Thus $\Gamma_1$ defined as (\ref{gama1}) is still the only heteroclinic orbit from $\mathbf{q}$ to $\Upsilon$.
\\\indent a) If $q_3=d-\sqrt{\rho}$ and $\mathbf{p}_0\in [\mathbf{x}_-,\mathbf{x}_+).$ \\
\indent Then (\ref{p0}) holds still.
Let $\Gamma_2$ be defined as (\ref{gama2}). Since the left subsystem in Theorem 2 is same as Theorem 1, the verification for $\{\phi_A(t,\mathbf{p}_0)|-\infty<t<0\}\subset \Sigma_-$
is same as (\ref{M6})  in Theorem 1. Now we will prove
\begin{eqnarray}\label{M15}
\{\phi_B(t,\mathbf{p}_0)|0<t<+\infty\} \subset \Sigma^+.
 \end{eqnarray}
Obviously, $\{\phi_B(t, \mathbf{p}_0)|t>0\}$ is contained in  $W^s(\mathbf{q})$.
From $B=\left(
 \begin{array}{ccc}
    b_{11} &b_{12} &0\\
    b_{21} & b_{22}& 0\\
    0 & 0& \lambda
  \end{array}\right)$ and condition (H2), we know, in $W^s(\mathbf{q})$, the flow of (\ref{M3}) is determined only by the following planar system

   $$\left(\begin{array}{c}
   \dot{x}_1\\
    \dot{x}_2
  \end{array}\right)=B_0
  \left(\begin{array}{c}x_1-q_1\\x_2-q_2\end{array}\right),$$
where $B_0=\left(\begin{array}{cc}b_{11}&b_{12}\\b_{21}&b_{22}\end{array}\right)$ with its eigenvalues being $\alpha\pm \beta i$ ($\alpha<0,\beta>0$).
From Lemma 3 and the definitions of $\mathbf{x}_-$ (see (\ref{Sigma5})) and $\mathbf{x}_+$(see Theorem 2), it is readily obtained that (\ref{M15}) holds
by $\mathbf{p}_0\in [\mathbf{x}_-,\mathbf{x}_+)$. Thus, $\Gamma_2$ is the only heteroclinic orbit from $\Upsilon$ to $\mathbf{q}$.
Therefore, $$\Upsilon\cup \Gamma_2 \cup \mathbf{q} \cup \Gamma_1$$ is the only heteroclinic cycle connecting $\Upsilon$ and $\mathbf{q}$.\\\\\
\indent b) If $q_3=d+\sqrt{\rho}$, $\omega^2\rho<\mu^2(d^2-\rho)$ and $\mathbf{p}_1\in [\mathbf{x}_-,\mathbf{x}_+)$.\\
\indent In this subcase,$$W^u(\Upsilon)\cap \Sigma \cap W^s(\mathbf{q})=\{\mathbf{p}_1\}.$$
Moreover, the proof for $\{\phi_A(t,\mathbf{p}_1)|-\infty<t<0\}\subset \Sigma_-$ is the same as the proof of (\ref{M11}),
 the  proof for $\{\phi_B(t,\mathbf{p}_1)|0<t<+\infty\}\subset \Sigma^+$ is similar to the proof of (\ref{M15}). We omit them here. Then,
the only heteroclinic cycle connecting $\Upsilon$ and $\mathbf{q}$ in this subcase can be expressed by
$$\Upsilon\cup \Gamma_3 \cup \mathbf{q} \cup \Gamma_1,$$
where $\Gamma_1$ and $\Gamma_3$ are defined as (\ref{gama1}) and (\ref{gama3}).\\
\indent c) if  $d-\sqrt{\rho}<q_3<d+\sqrt{\rho}$, $\omega^2\rho<\mu^2(d^2-\rho)$ and $\mathbf{p}_\pm\in [\mathbf{x}_-,\mathbf{x}_+)$.\\
\indent In this subcase, we have
$$W^u(\Upsilon)\cap \Sigma \cap W^s(q)=\{\mathbf{p}_+,\mathbf{p}_-\}.$$
The proofs for $\{\phi_A(t,\mathbf{p}_\pm)|t<0\}\subset \Sigma_-$ and $\{\phi_B(t,\mathbf{p}_\pm)|t>0\}\subset \Sigma^+$ is similar to the proofs of (\ref{M11}) and (\ref{M15}), respectively. We omit them here.
 Then,
the only two heteroclinic cycles connecting $\Upsilon$ and $\mathbf{q}$ in this subcase can be expressed by
$$\Upsilon\cup \Gamma_\pm \cup \mathbf{q} \cup \Gamma_1,$$
where $\Gamma_1$ and $\Gamma_\pm$ are defined as (\ref{gama1}) and (\ref{Gammapn}), respectively.~~~~~$\square$\vspace{3mm}\\
\textbf{The proof of (ii)( When $0<d^2-\rho<\frac{\omega^2}{4d^2}$).}\vspace{3mm}\\
\indent Combining the proof of (ii) of Theorem 1 and the proof of (i) of Theorem 2, this proof is easy to achieve here. We omit it for simplification. $\square$.

\section{Examples}
\begin{figure}[!htbp]
\renewcommand{\figurename}{Fig.}
 \begin{minipage}{0.25\textwidth}
\caption{A heteroclinic cycle connecting $\Upsilon$ and $\mathbf{q}$ in Example 1. }\label{fig3}
\end{minipage}
\hspace{0cm}
 \begin{minipage}{0.75\textwidth}
\includegraphics[width=9cm]{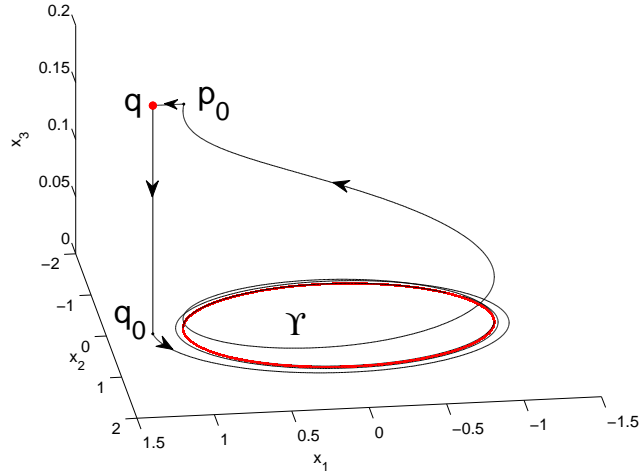}
\end{minipage}
\end{figure}

\noindent \textbf{Example 1:Heteroclinic cycle connecting a pure saddle (i.e., all eigenvalues are real) and a saddle periodic orbit}. \\
For system (\ref{M1}), let $$
 A=\left(
 \begin{array}{ccc}
    \rho &-\omega &0\\
    \omega & \rho& 0\\
    0 & 0& \mu
  \end{array}\right)=\left(
 \begin{array}{ccc}
    1 &-10 &0\\
    10 & 1& 0\\
    0 & 0& 5
  \end{array}\right),B=\left(
 \begin{array}{ccc}
    b_{11} &b_{12} &0\\
    b_{21} & b_{22}& 0\\
    0 & 0& \lambda
  \end{array}\right)=\left(
 \begin{array}{ccc}
   -2 &1&0\\
    0 & -1& 0\\
    0 & 0& 2
  \end{array}\right),$$
$$~\mathbf{q}=(q_1,q_2,q_3)^T=(1.2,0,0.2)^T,~d=1.2.$$
  Then, \begin{eqnarray*}
\mathbf{p}_0=(\sqrt{\rho},0,d-\sqrt{\rho})^T=(1,0,0.2)^T
, \mathbf{q}_0=(d,q_2,0)^T=(1.2,0,0)^T.
\end{eqnarray*}
 $$\sigma_+=\frac{-\omega+\sqrt{\omega^2-4d^2(d^2-\rho)}}{2d}=-0.05314,~~\mathbf{\mathbf{v}}_1=(d,,0)^T=(1.2,-0.05314,0)^T.$$
 $$\rho=1<1.44=d^2, \textbf{c}^T\mathbf{q}=1.4>1.2=d,q_1=d=1.2.$$
Moreover, the eigenvalues of $B$: $\lambda_1=-2<0,\lambda_2=-1<0$ and $\lambda=2>0$, thus, $B$ satisfies (H1) and (H3). And
$$0<d^2-\rho=0.44<\frac{625}{36}=\frac{\omega^2}{4d^2}.$$
In addition,
$\mathbf{v}^*=(d,v^*_2,0)^T=(1.2,2.363,0)^T$ by numerical calculation.
 $$v^*_2=2.363>-0.05314=\sigma_+,~~\sigma_+< q_2=0< v^*_2.$$
 $$q_3=d-\sqrt{\rho}=0.2,~~ \mathbf{c}^TB(\mathbf{p}_0-\mathbf{q})=0.4\geq0.$$
 According to conclusion a) in case (ii) in Theorem 1, there exists only one heteroclinic cycle connecting periodic orbit $\Upsilon=\{(x,y,z)\in \mathbb{R}^3|x^2+y^2=1,z=0\}$
 and equilibrium $\mathbf{q}$ as shown in Figure \ref{fig3}.
\\ \noindent \textbf{Example 2: Heteroclinic cycle connecting a saddle-focus and a saddle periodic orbit}.\\
 For system (\ref{M1}), let \begin{eqnarray*}
 A=\left(
 \begin{array}{ccc}
    1 &-35^{\frac{1}{2}} &0\\
    35^{\frac{1}{2}} & 1& 0\\
    0 & 0& 5
  \end{array}\right),B=\left(
 \begin{array}{ccc}
   -0.5 &4 &0\\
    -4 & -0.5& 0\\
    0 & 0& 2
  \end{array}\right),
  \end{eqnarray*}
 $$d=(35/11)^{\frac{1}{2}}, ~\mathbf{q}=(d,-4.5,d+1)^T$$
 Then $$\mathbf{p}_1=(-\sqrt{\rho},0,d+\sqrt{\rho})^T=(-1,0,d+1),$$
  $$\sigma_+=\frac{-\omega+\sqrt{\omega^2-4d^2(d^2-\rho)}}{2d}=-0.9045,\sigma_-=\frac{-\omega-\sqrt{\omega^2-4d^2(d^2-\rho)}}{2d}=-2.4121.$$
  \begin{figure}[!htbp]
\renewcommand{\figurename}{Fig.}
 \begin{minipage}{0.25\textwidth}
\caption{A heteroclinic cycle connecting periodic orbit $\Upsilon=\{(x_1,x_2,x_3)^T|x_1^2+x_2^2=1,x_3=0\}$ and saddle-focus $\mathbf{q}$ in Example 2. }\label{fig4}
\end{minipage}
\hspace{0cm}
 \begin{minipage}{0.75\textwidth}
\includegraphics[width=10cm]{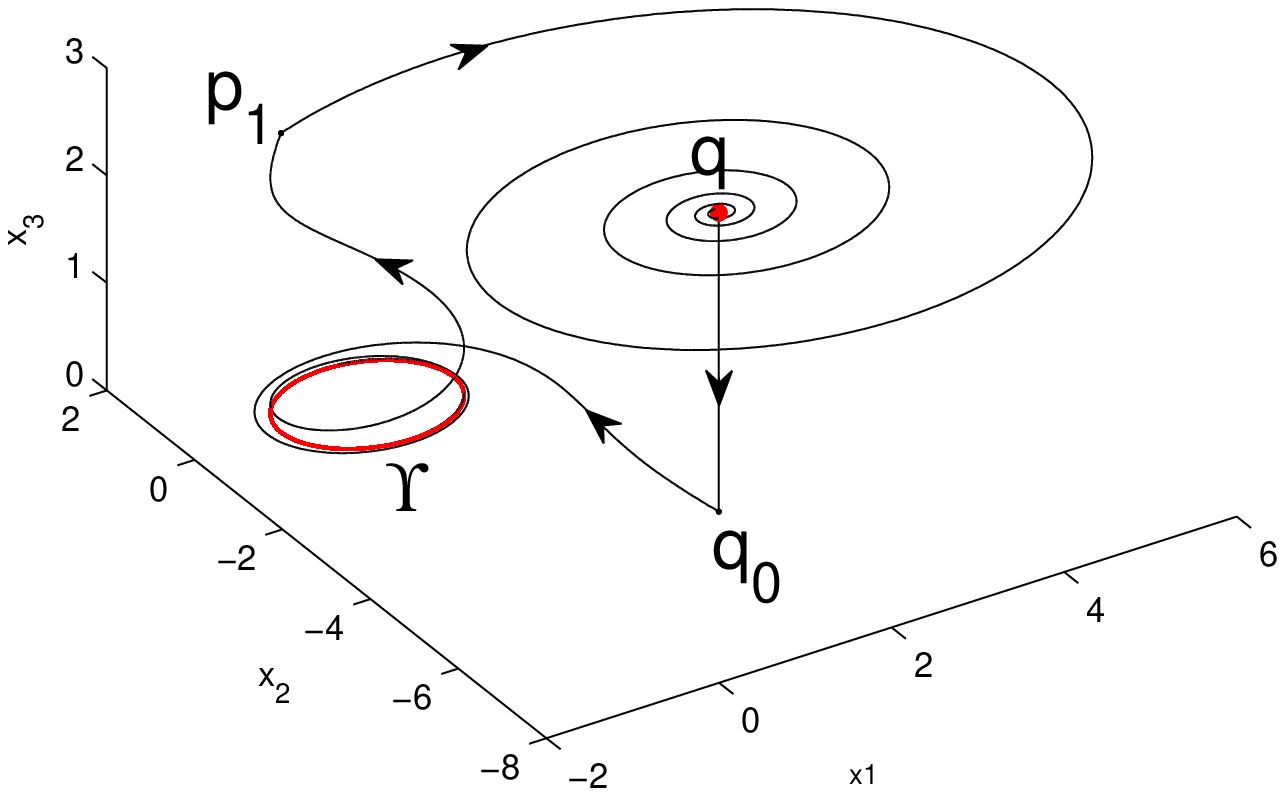}
\end{minipage}
\end{figure}
$$\rho=1<\frac{35}{11}=d^2,~~\textbf{c}^T\mathbf{q}=2d+1>d,~~q_1=d=(35/11)^{\frac{1}{2}},$$
$$\mathbf{\mathbf{v}}_1=(d,\sigma_+,0)^T=(d, -0.9045,0)^T,$$
Then, the eigenvalues of $B$ are $-0.5\pm 4i$ and 2. Thus, $B$ satisfies (H2) and (H3). And $$0<d^2-\rho=\frac{24}{11}<\frac{11}{4}=\frac{\omega^2}{4d^2}.$$
By numerical calculation, $\mathbf{v}^*=(d,v^*_2,0)^T=(d, -4.4162,0)^T$. Then,
 $$v^*_2=-4.4162<-2.4124=\sigma_-,~~q_2=-4.5<v^*_2=-4.4162.$$
In addition, $$q_3=d+\sqrt{\rho}=d+1,~ \omega^2\rho=35<\frac{600}{11}=\mu^2(d^2-\rho).$$ Meanwhile,  $\mathbf{x}_-=\frac{d-\mathbf{c}^T\mathbf{q}}{\mathbf{c}^TB^{-1}\mathbf{c}^\bot} B^{-1}\mathbf{c}^\bot+\mathbf{q }=(-1, -4.848,d+1)^T$ and $\mathbf{x}^*=(-1,0.0476,d+1)^T$ by numerical calculation.
Thus $$\mathbf{p}_1=(-\sqrt{\rho},0,d+\sqrt{\rho})^T=(-1,0,d+1)\in [\mathbf{x}_-,\mathbf{x}^*).$$
According to conclusion b) in case (ii) of Theorem 2, there exists only one heteroclinic cycle connecting periodic orbit $\Upsilon=\{(x,y,z)\in \mathbb{R}^3|x^2+y^2=\rho=1,z=0\}$
 and equilibrium $\mathbf{q}$ as shown in Figure \ref{fig4}.

  \begin{figure}[!htbp]
\renewcommand{\figurename}{Fig.}
 \begin{minipage}{0.25\textwidth}
\caption{Two heteroclinic cycles each connecting periodic orbit $\Upsilon=\{(x_1,x_2,x_3)^T|x_1^2+x_2^2=1,x_3=0\}$ and saddle-focus $\mathbf{q}$ in Example 3. }\label{fig5}
\end{minipage}
\hspace{0cm}
 \begin{minipage}{0.75\textwidth}
\includegraphics[width=12cm]{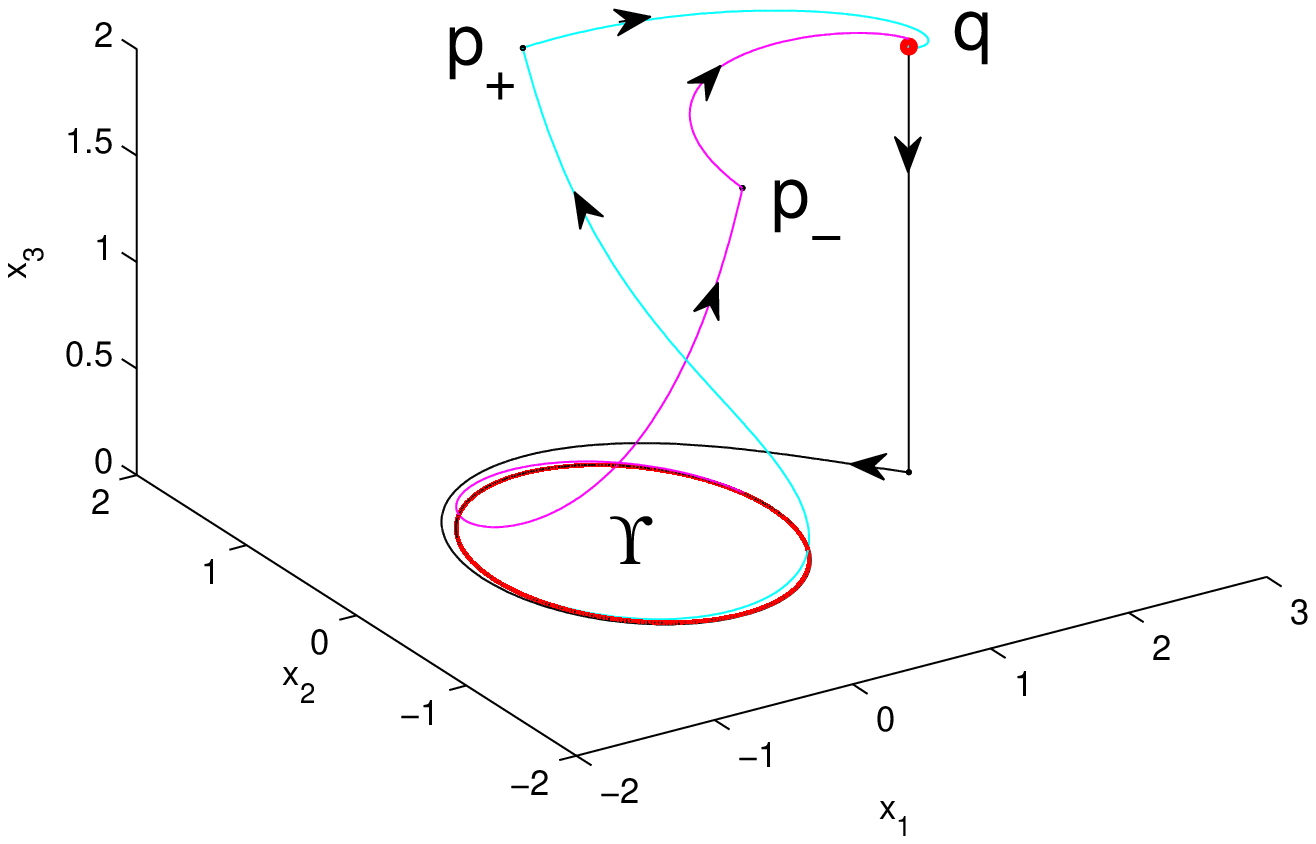}
\end{minipage}
\end{figure}
\noindent \textbf{Example 3: A pair of  heteroclinic cycles each connecting a saddle-focus and a saddle periodic orbits}.\\
 For system (1), let \begin{eqnarray*}
 A=\left(
 \begin{array}{ccc}
    1 &-3 &0\\
    3 & 1& 0\\
    0 & 0& 4
  \end{array}\right),B=\left(
 \begin{array}{ccc}
   -3.5 &6 &0\\
    -6 & -3.5& 0\\
    0 & 0& 2
  \end{array}\right),.
  \end{eqnarray*}
 $$d=2, ~\mathbf{q}=(2,0,2)^T.$$
  Then, $$\mathbf{p}_-=(0,-1,2)^T,~\mathbf{p}_+=(0,1,2)^T,~\mathbf{q}_0=(2,0,0)^T.$$
  $$\rho=1<4=d^2,~~\textbf{c}^T\mathbf{q}=4>2=d,~~q_1=d=2.$$
Moreover, the eigenvalues of $B$: $-3.5\pm 6i$, thus, $B$ satisfies (H2). And
 $$d^2-\rho=3>2.5=\frac{\omega^2}{4d^2}.$$
Furthermore,
\begin{center}
 $d-\sqrt{\rho}=1<q_3=2<3=d+\sqrt{\rho}$, $\omega^2\rho=36<48=\mu^2(d^2-\rho)$,
\end{center}
and  $\mathbf{x}_-=(0,-1.1667,2)^T$ and $\mathbf{x}^*=( 0,27.6586,2)^T$ by numerical calculation.
Thus $$\mathbf{p}_\pm\in (\mathbf{x}_-,\mathbf{x}^*).$$
 According to the conclusion c) in case (i) of Theorem 2, there exist only two heteroclinic cycles each connecting periodic orbit $\Upsilon=\{(x,y,z)\in \mathbb{R}^3|x^2+y^2=\rho=1,z=0\}$
 and saddle-focus $\mathbf{q}$ as shown in Figure \ref{fig5}.

\section{Conclusions}
In this paper, we have investigated the existence of  heteroclinic cycles connecting saddle periodic orbit and saddle equilibrium in a class of piecewise smooth systems. Two types of such heteroclinic cycles are constructed:
one is the heteroclinic cycle connecting a saddle point with only real eigenvalues and a saddle
periodic orbit, and the other is the heteroclinic cycle connecting a saddle-focus point and a saddle
periodic orbit. What's more, the main results are convenient and feasible to be used in the construction of piecewise smooth systems possessing such heteroclinic cycles, see these examples in Section 5.\\
\indent As mentioned in the Introduction section, for flows, the singular cycles connecting saddle periodic orbit and saddle equilibrium can potentially result in the so-called singular horseshoe, which implies the existence of a non-uniformly hyperbolic chaotic invariant set.  Maybe, It can be conjectured that the singular cycles studied in this paper for the piecewise smooth systems may also result in the emerge of the singular horseshoe under some conditions. We will investigate this interesting issue in the future work.
\section*{Acknowledgments}
\noindent The first author takes this opportunity to thank Prof. Sebastian Walther (Department of Mathematics, RWTH Aachen University, Germany) for providing a good research environment and office space during the visit in the RWTH Aachen University from July 2018 to July 2019.\\ This study was supported by the National Natural Science Foundation of China (Grant numbers 11702077 and 11472111), the Natural Science Foundation of Anhui Province (Grant number 1708085QA12), the Overseas Visiting and Training Foundation of Outstanding Young Talents in the Universities of Anhui Province (Grant number gxgwfx2018070), the Fostering Master's Degree Empowerment Point Project of Hefei University(Grant number 2018xs03) and the Anhui Provincial Teaching Studio (Grant number 2014MSGZS159).

\section*{References}
\bibliographystyle{elsarticle-num}
\bibliography{refer_Van}

\end{document}